\documentclass{amsart}

\newtheorem{theorem}{Theorem}[section]
\newtheorem{lemma}[theorem]{Lemma}
\newtheorem{corollary}[theorem]{Corollary}
\newtheorem{proposition}[theorem]{Proposition}

\theoremstyle{remark}
\newtheorem{rem}{Remark}[section]
\numberwithin{equation}{section}

\def\mod#1{{\ifmmode\text{\rm\ (mod~$#1$)}
\else\discretionary{}{}{\hbox{ }}\rm(mod~$#1$)\fi}}

\begin{document}

\title[Cubic Thue equations]
{Extremal families of cubic Thue equations}

\author{Michael A. Bennett and Amir Ghadermarzi}
\address{Department of Mathematics, University of British Columbia}

\email{bennett@math.ubc.ca}
\email{amir@math.ubc.ca}
\thanks{The first author was supported in part by NSERC}
\dedicatory{On the occasion of Axel Thue's 150th birthday}
\subjclass{Primary 11D25, Secondary 11E76}

\date{}

\keywords{Thue equations, binary cubic forms}

\begin{abstract}
We exactly determine the integral solutions to a previously untreated infinite family of  cubic Thue equations of the form $F(x,y)=1$ with at least $5$ such solutions. Our approach combines elementary arguments, with lower bounds for linear forms in logarithms and lattice-basis reduction.
\end{abstract}

\maketitle

\section{Introduction}

If $F(x,y)$ is an irreducible binary form (i.e. homogenous polynomial) with integer coefficients and degree $n \geq 3$, and $m$ is a nonzero integer, the Diophantine equation $F(x,y)=m$ is called a {\it Thue equation} in honour of Axel Thue \cite{Thu} who proved, more than a century ago, that the associated number of solutions in integers $x$ and $y$ is finite.  Such equations arise somewhat naturally in a wide variety of number theoretic contexts, including questions about the existence of primitive divisors in binary recurrence sequences and in the only known algorithm (in the strict sense of the term) for finding integral points on genus one curves over $\mathbb{Q}$. Our current understanding of such equations is a rather refined one and we now have excellent  upper bounds upon their number of solutions, depending only on $n$ and $m$ (and, in particular, not upon the coefficients of the form; see Evertse \cite{Eve} and Bombieri and Schmidt \cite{BS}). These bounds arise from treating the apparently special case $m=1$ and applying ``lifting'' arguments, dating back to Lagrange (see page 673 of Dickson \cite{Dic} and also Mahler \cite{Mah}). 

The case when $n=3$, i.e. that of cubic forms, is a relatively simple one and is well understood. This is primarily because  cubic Thue equations may be attacked with a wide variety of approaches that do not apparently generalize to those of higher degree. The particular case of cubic Thue equations of the shape $F(x,y)=1$ where the form $F$ has a negative discriminant is especially simple, as the presence of a single fundamental unit in the corresponding cubic field makes application of Skolem's $p$-adic method or similar local techniques relatively routine. Using such an approach (where we write $N_F$ for the number of integral solutions to the equation $F(x,y)=1$ and $D_F$ for the discriminant of $F$), Delone \cite{Del} and Nagell \cite{Nag}
 independently  proved

\begin{theorem} \label{DelNag} (Delone and Nagell)
If $F$ is an irreducible binary cubic form with integer coefficients and
$D_F < 0$, then $N_F \leq 5$. Moreover, if $N_F = 5$, then $F$ is $GL_2 (\mathbb{Z})$-equivalent to
$$
x^3-xy^2+y^3,
$$
with $D_F = -23$
and, if $N_F=4$, then $F$ is $GL_2 (\mathbb{Z})$-equivalent to either
$$
x^3+xy^2+y^3 \mbox{ or }  x^3-x^2y+xy^2+y^3,
$$
with discriminant $-31$ or $-44$, respectively.
\end{theorem}

In the case of cubic forms of positive discriminant, the situation is rather more complicated and there are a number of forms for which $N_F$ exceeds $5$.
The following table collects representatives of all known equivalence classes of irreducible cubic forms
for which $N_F \geq 6$.

$$
\begin{array}{|cccl|} \hline
F(x,y) & D_F & N_F & \mbox{References} \\ \hline
x^3-x^2y-2xy^2+y^3 & 49 & 9 & \mbox{\cite{Bau}, \cite{GaS}, \cite{Lju}, \cite{PeS}} \\
x^3-3xy^2+y^3 & 81 & 6 & \mbox{\cite{GaS}, \cite{Lju}, \cite{Tza}} \\
x^3-4xy^2+y^3 & 229 & 6 & \mbox{\cite{Bre}, \cite{GaS}, \cite{PeS}} \\
x^3-5xy^2+3y^3 & 257 & 6 & \mbox{\cite{GaS}} \\
x^3+2x^2y-5xy^2+y^3 & 361 & 6 & \mbox{\cite{GaS}} \\ \hline
\end{array}
$$
\vskip1ex \noindent
Presumably, we always have $N_F \leq 9$. A striking theorem of Okazaki \cite{Oka} is that $N_F \leq 7$ provided $D_F$ is suitable large (see also Akhtari \cite{Akh}); we have $N_F \leq 10$ in all cases, via an old result of the first author \cite{Ben2}.

A stronger conjecture, due originally to Nagell \cite{Nag2}
and subsequently refined by Peth\H{o} \cite{Pet} and Lippok \cite{Lip},  is that the forms
in the above table are, up to equivalence, the only irreducible cubics with $N_F \geq 6$, so that, for 
all other classes, we have $N_F \leq 5$. If true, this
upper bound for $N_F$  is sharp, as we know of a number of infinite families of cubic forms where the number of integer solutions to the corresponding equation $F(x,y)=1$ is at least $5$. Let us define
\begin{equation} \label{five-1}
F_{1,t}(x,y)=x^3-(t+1)x^2y+txy^2+y^3,
\end{equation}
\begin{equation} \label{five-2}
F_{2,t}(x,y)=x^3-t^2 xy^2+y^3
\end{equation}
and
\begin{equation} \label{five-3}
F_{3,t}(x,y) = x^3 - (t^4 -t)x^2y + (t^5-2t^2)xy^2 + y^3.
\end{equation}
In each case, $F_i(x,y)$ is irreducible  over $\mathbb{Q} [x,y]$, with $N_{F_{i,t}} \geq 5$, at least provided we exclude ``small'' values of the integer parameter $t$ (the assumption that $|t| > 2$ is sufficient). 

It is plausible to believe that (\ref{five-1}), (\ref{five-2}) and (\ref{five-3})  represent the only infinite cubic families for which the corresponding Thue equation $F(x,y)=1$ has five or more integral solutions, partially explaining our (misleading)  title. The evidence does not admittedly seem especially compelling, one way or the other. That being said, we will take this opportunity to  solve those equations corresponding to the third family, which has not previously been treated in the literature. We note that the equation $F_{1,t}(x,y)=1$ was completely solved by Lee \cite{Lee}, Mignotte and Tzanakis \cite{MiT} and Mignotte \cite{Mig3}, while $F_{2,t}(x,y)=1$ has been treated for $|t| \geq 1.35 \cdot 10^{14}$ by Wakabayashi \cite{Wa5}. We prove

\begin{theorem} \label{main-thm}
If $t$ is an integer, then the Diophantine equation
\begin{equation} \label{main-eq}
F_{3,t}(x,y) = x^3 - (t^4 -t)x^2y + (t^5-2t^2)xy^2 + y^3 = 1,
\end{equation}
has only the integer solutions  
$$
(x,y) \in \{  (1,0), (0,1), (t,1),  (t^4-2t,1),  \left( 1-t^3,t^8-3t^5+3t^2 \right) \},
$$
unless $t=-1$, in which case there is an additional solution given by $(x,y)=(6,-5)$.
\end{theorem}

Our argument follows the now-traditional approach originated by Thomas \cite{Tho}, \cite{Th93}. A nice survey of families of Thue equations solved to date (by these and other methods) can be found in Heuberger \cite{Heu} (see also \cite{Heu2} and \cite{HTZ} for other good expositions along these lines).

\section{The equation $F_{3,t}(x,y)=1$ : units in cubic fields}

For  the remainder of the paper, we will concern ourselves with the parametric  equation (\ref{main-eq}).
Note that $N_{F_{3,t}} \geq 5$ (with solutions given in the statement of Theorem \ref{main-thm}) for $t \not\in \{ 0, 1  \}$.
We have
$$
D_{F_{3,t}} = t^{18}-10t^{15}+41t^{12}-90t^9+102t^6-40t^3-27,
$$
which is positive, except for $t \in \{  0, 1 \}$. That $N_{F_{3,t}} \geq 5$ was noted by Ziegler \cite{Zie} who observed that one also has $N_{F_{4,t}} \geq 5$, for the families of forms given by 
$$
F_{4,t}(x,y)=x^3-(t^4+4t)x^2y+(t^5+3t^2)xy^2+y^3.
$$
We have $F_{3,-t}(x,y) \sim F_{4,t}(x,y)$  under $GL_2(\mathbb{Z})$ action, since  $F_{3,-t}(x-ty,y)=F_{4,t}(x,y)$.

Let us suppose for the next few sections that $t \geq 10$; our argument for negative values of $t$ is very similar and will not be reproduced here, while the treatment of  ``small'' values of $t$ requires rather different techniques. 
Defining $P(x)=F_{3,t}(x,1)$,
then $P(x)$ has three real roots, which we denote by $\theta_1 < \theta_2 < \theta_3$. By studying the sign of $P(x)$, we can deduce the following expressions for these roots :
\begin{equation} \label{estimation}
 \theta_1 = - \frac{1}{t^5}-\frac{2}{t^8}-\frac{\kappa_1}{t^{11}}, \;  \;
 \theta_2 = t+\frac{1}{t^5}+ \frac{3}{t^8}+\frac{\kappa_2}{t^{11}}  \; \mbox{ and } \; 
\theta_3 = t^4-2t-\frac{1}{t^8} -\frac{\kappa_3}{t^{11}}. 
\end{equation}
Here, the $\kappa_i$ are certain real numbers with $\kappa_1 \in (3,3.1)$, $\kappa_2 \in (8,8.03)$  and $\kappa_3 \in (5,5.02)$.

Suppose that $(x,y)$ is a solution to equation (\ref{main-eq}). It follows, for each $i \in \{ 1, 2, 3\}$, that $x- y \theta_i$ is a unit in $\mathbb{Q}(\theta_i)$  
and hence, since Lemma 4.11 of \cite{Th93} implies that $t-\theta_i$ and $\theta_i$ form a pair of fundamental units in $\mathbb{Q}(\theta_i)$, we may write 
\begin{equation}\label{unit}
{x - y \theta_i = (-1)^\delta (t - \theta_i)^n \theta_i^{-m} \ \ \text{for} \ m,n \in \mathbb{Z}}, \; \delta \in \{ 0, 1 \}.
\end{equation}
In particular, we have
 \begin{equation} \label{good}
  \frac{x-y\theta_3}{x-y\theta_2}=\left(\frac{t-\theta_3}{t-\theta_2} \right)^n \left(\frac{\theta_3}{\theta_2}\right)^{-m} > 0,
  \end{equation}
and so, from the fact that
$$
  (x-y\theta_1)(x-y\theta_2) (x-y\theta_3)=1, 
$$
  we conclude that $x-y\theta_1 >0$ (whereby at least one of $x$ and $y$ is positive). 

If $|y| \leq 1$, we observe after a little work, since we assume $t \geq 10$, that
$$
(x,y) \in \{  (1,0), (0,1), (t,1),  (t^4-2t,1) \}.
$$
Let us therefore suppose that $(x,y)$ is a solution in integers to equation (\ref{main-eq}), with $|y| \geq 2$.
A routine calculation ensures, provided $t \geq 10$, that necessarily $x/y$ lies in one of the intervals
$$
I_1=\left( - \frac{1.13}{t^5}, \frac{-1 + \frac{1}{|y|^3}}{t^5}  \right), \; 
I_2=\left( t+\frac{1 - \frac{1}{|y|^3}}{t^5}, t + \frac{1.13}{t^5} \right)
$$
or
$$
I_3 = \left( t^4-2 t - \frac{1.13}{t^8} , t^4 - 2 t - \frac{1 - \frac{1}{|y|^3}}{t^8} \right).
$$
As is readily observed, these intervals are disjoint and hence we will call solutions $(x,y)$ to equation (\ref{main-eq}) with $x/y$ in $I_1$, $I_2$ and $I_3$ solutions of types I, II and III, respectively. It is easy to show (and valuable to note, for later use) that $|x-y \theta_i|<1$ for $i=1, 2$ or $3$, provided $(x,y)$ is a solution of type I, II or III, respectively.
  
  From here, we will proceed as follows. We first will use elementary arguments and a careful analysis of equation (\ref{good}) to deduce lower bounds upon $\max \{ |m|, |n| \}$ in (\ref{unit}). We then combine this information with lower bounds upon
  \begin{equation} \label{rip}
  \left|  \frac{x-y\theta_i}{x-y\theta_j} \right|
  \end{equation}
  for suitably chosen pairs $(i,j)$, depending on the solution type of $(x,y)$. The latter bounds arise from invoking lower bounds for linear forms in logarithms of algebraic numbers \`a la Baker and yield upper bounds for $t$. Finally, we appeal to a classical lemma of Baker and Davenport  \cite{BaDa} (in essence, a simple version of Lenstra-Lenstra-Lovacz lattice basis reduction) to treat (most of)  the remaining values of $t$.
   
\section{Upper bounds for $\max \{ |m|, |n| \}$  in  equation (\ref{good})}

In general, it is always possible to reduce the problem of solving families of Thue equations to that of treating unit equations similar to (\ref{unit}). A crucial step in solving such families  is to ensure suitably rapid growth (Thomas \cite{Th93} terms this {\it stable growth}) of the exponents of the fundamental units. That is the content of this section. Our arguments vary somewhat depending on the solution type of $(x,y)$; we will treat each in turn.

  \subsection{Solutions of type I}
  Let us suppose first that $x/y \in I_1$, i.e. that $x/y$ satisfies
  \begin{equation} \label{cool-daddy1}
  - \frac{1.13}{t^5} < \frac{x}{y} <  \frac{-1 + \frac{1}{|y|^3}}{t^5}.
   \end{equation}
If we have $y<0$ and $x>0$ then from 
$x-y\theta_1 >0$, it follows that $x/y \in ( - \frac{1.13}{t^5}, \theta_1)$.
   From (\ref{estimation}), we thus have
$$
\frac{t^3-2+\frac{1}{t^6} + \frac{2-\kappa_3}{t^9}+\frac{\kappa_1}{t^{12}}}{1+\frac{2.13}{t^6}+\frac{\kappa_2}{t^9}} <  \frac{x-y\theta_3}{x-y\theta_2} < \frac{t^3-2+\frac{1.13}{t^6} - \frac{\kappa_3}{t^9}}{1+\frac{2}{t^6}+\frac{2+\kappa_2}{t^9}+\frac{\kappa_1}{t^{12}}}.
 $$
 If, on the other hand, we have $y>0$ (so that $y \geq 2$) and $x < 0$, then $x/y \in\left( \theta_1, \frac{-0.875}{t^5} \right)$, whereby
 $$
\frac{t^3-2+\frac{0.875}{t^6}- \frac{\kappa_3}{t^9}}{1+\frac{2}{t^6}+\frac{2+\kappa_2}{t^9}+\frac{\kappa_1}{t^{12}}} <  \frac{x-y\theta_3}{x-y\theta_2} < \frac{t^3-2+\frac{1}{t^6} + \frac{2-\kappa_3}{t^9}+\frac{\kappa_1}{t^{12}}}{1+\frac{1.875}{t^6}+\frac{\kappa_2}{t^9}}.
 $$
In either case, we thus may write
 $$
 \frac{x-y\theta_3}{x-y\theta_2} = t^3-2 - \frac{\kappa_4}{t^3},
 $$
 where $\kappa_4 \in (1.8,2.2)$. Arguing similarly, we find that
 $$
 \log \left(  \frac{t-\theta_3}{t-\theta_2}  \right) = 9 \log t - \frac{6}{t^3} - \frac{\kappa_5}{t^6},
 $$
 $$
 \log \left( \theta_3/\theta_2 \right) = 3 \log t  - \frac{2}{t^3}  - \frac{\kappa_6}{t^6} 
 $$
 and
 $$
 \log \left(  \frac{x-y\theta_3}{x-y\theta_2}  \right) = 3 \log t  - \frac{2}{t^3}  - \frac{\kappa_7}{t^6},
 $$
 where $ \kappa_5 \in (7.99,8.03)$, $\kappa_6 \in (3,3.01)$ and  $\kappa_7 \in (3.8,4.3)$.
 
 We may thus conclude from (\ref{good}) that
 \begin{equation} \label{help1}
k:=3n-m-1 = \frac{1}{3 \log t} \left( \frac{6 n-2m-2}{t^3} 
+ \frac{ \kappa_5 n  - \kappa_6m - \kappa_7}{t^6} \right).
\end{equation}
If $k=0$ (so that $m=3n-1$), then it follows that
$$
0=\kappa_5 n - \kappa_6 (3n-1) - \kappa_7 = (\kappa_5- 3 \kappa_6) n + \kappa_6 - \kappa_7
$$
and so 
$$
n = \frac{\kappa_6 - \kappa_7}{3 \kappa_6 - \kappa_5} \in (-1.25,-0.8), \; \mbox{ i.e. } n=-1, \; m=-4.
$$
Equation (\ref{unit}) thus implies, after a little work, that
$$
(x,y) = \left( 1-t^3,t^8-3t^5+3t^2 \right).
$$
 We may thus suppose that $k$ is a nonzero integer whence, from equation (\ref{help1}), we have
 $$
k \left( 3 t^6 \log t - 2 t^3 \right) =  \kappa_5 n  - \kappa_6 m- \kappa_7 = 
\frac{\kappa_5}{3} k + \frac{\kappa_5 - 3 \kappa_6}{3} m + \frac{\kappa_5 - 3 \kappa_7}{3} .
 $$
  It follows that 
  \begin{equation} \label{emm-1}
  m = \frac{k}{\kappa_5 - 3 \kappa_6 } \left( 9t^6 \log t - 6 t^3 - \kappa_5 \right) 
  +\frac{3 \kappa_7-\kappa_5}{\kappa_5 - 3 \kappa_6}.
  \end{equation}
  We may conclude from this that $m$ and $k$ are necessarily of opposite sign and, since $t \geq 10$,
  that
  $$
  |m| > 8.6 \, |k| \, t^6 \log t.
  $$
  From $k=3n-m-1$, $m$ and $n$ thus have the same sign whereby, considering (\ref{unit}) and the fact that $|x - y \theta_1|<1$, we may conclude that $m < n < 0$, $|m| > 3 |n|$ and, since $k$ is nonzero,  
  \begin{equation} \label{emm-2}
  |m| = \max \{ |m|, |n| \} > 8.6 \, t^6 \log t.
  \end{equation}

  \subsection{Solutions of type II}
  Let us next suppose that $x/y$ satisfies
  \begin{equation} \label{cool-daddy}
   t+\frac{1}{t^6} < \frac{x}{y} < t+\frac{2}{t^5}.
   \end{equation}
   Since at least one of $x$ and $y$ is positive, it follows that both are positive.
  We therefore have (again appealing to (\ref{estimation}) and (\ref{unit}))
   $$
   x-y\theta_1 > 0, \; \; x-y\theta_2 < 0  \; \mbox{ and } \; x-y\theta_3 <0.
   $$
 It follows that 
     $$ 
     0 > \frac{x- y \theta_3}{x- y \theta_1}= \left(\frac{t-\theta_3}{t-\theta_1}\right)^n \left(\frac{\theta_3}{\theta_1} \right)^{-m}
     $$
    and hence $n$ and $m$ are of opposite parity.
  
   Using the inequalities in  (\ref{estimation}) and (\ref{cool-daddy}), we find that
    $$
  \frac{t^3-3-\frac{2}{t^6} -\frac{\kappa_3}{t^9}}{1+
    \frac{3}{t^6}+\frac{2}{t^9}+\frac{\kappa_1}{t^{12}}} <  \frac{y \theta_3 -x}{x- y \theta_1} < \frac{t^3-3-\frac{1}{t^7} -\frac{\kappa_3}{t^9}}{1+
    \frac{1}{t^6}+\frac{1}{t^7}+\frac{2}{t^9}+\frac{\kappa_1}{t^{12}}}
    $$
    and hence
    $$
    \frac{y \theta_3 -x}{x- y \theta_1} = t^3-3-\frac{\kappa_8}{t^3}, \; \mbox{ where } \kappa_8 \in (0, 3.1).
    $$
    Similarly, 
$$
\frac{\theta_3-t}{t-\theta_1} = t^3-3-\frac{\kappa_9}{t^3}, \; \mbox{ where } \kappa_9 \in (0, 1.1) 
$$
and
$$
 \frac{\theta_3}{\lvert \theta_1 \rvert } = t^9-4t^6+\kappa_{10} t^3 \mbox{ with } \kappa_{10} \in (4.9,5).
$$
 We thus have
  \begin{equation} \label{model}
    \log \left(t^3-3-\frac{\kappa_8}{t^3} \right) = n \log \left( t^3-3-\frac{\kappa_9}{t^3} \right ) -m \log \left(t^9-4t^6+\kappa_{10} t^3 \right ).
  \end{equation}   
 If $m=0$, it thus follows that  $n=1$ which, with positive sign, leads us to the solution $(x,y)=(t,1)$.
 If $m<0$ then necessarily $n<0$ and  so 
 $$
  \lvert x-y\theta_2 \rvert = \lvert (t-\theta_2)^n \theta_2^{-m}\rvert >1,
  $$
  contradicting 
 (\ref{main-eq}). We may thus conclude that $m>0$ (whereby also $n > 0$).
 Since
 $$
 \log \left(t^3-3-\frac{\kappa_8}{t^3} \right) = 3 \log t - \frac{3}{t^3} - \frac{\kappa_{11}}{t^6} \mbox{ for } \kappa_{11} \in (4.5,7.7),
 $$
 $$
 \log \left(t^3-3-\frac{\kappa_9}{t^3} \right) = 3 \log t - \frac{3}{t^3} - \frac{\kappa_{12}}{t^6} \mbox{ for } \kappa_{12} \in (4.5,5.7)
 $$
 and
 $$
 \log \left(t^9-4t^6+\kappa_{10} t^3 \right) = 9 \log t - \frac{4}{t^3} - \frac{\kappa_{13}}{t^6} \mbox{ for } \kappa_{13} \in (2.9,3.1),
 $$
we may conclude that
\begin{equation} \label{help}
k := n-3m-1 = \frac{1}{3 \log t} \left( \frac{3n-4m-3}{t^3} 
+\frac{ \kappa_{12} n - \kappa_{13} m - \kappa_{11}}{t^6} \right).
\end{equation}
 Since equation (\ref{model}) readily implies that $n > 2m$, (\ref{help}) and the fact that $n$ and $m$ are of opposite parity allows us to conclude that $k$ is an even positive integer (so that, in particular,  $n \geq 3m+3$).
 Arguing crudely, (\ref{help}) thus implies that 
 $$
 2 \leq \frac{1}{3 \log t } \left( \frac{5n}{3 t^3} + \frac{5n}{t^6} \right)
 $$
 and hence
 \begin{equation} \label{lower}
 n = \max \{ |m|, |n| \}  > 3.5 \, t^3 \log t.
 \end{equation}
 
\subsection{Solutions of type III}

Suppose now that we have a solution $(x,y)$ with $x/y$ in $I_3$, i.e.
$$
 t^4-2 t - \frac{1.13}{t^8}  < \frac{x}{y} < t^4 - 2 t - \frac{1 - \frac{1}{|y|^3}}{t^8}.
$$
 As with solutions of type II, we may suppose that both $x$ and $y$ are positive integers (with, say, $y \geq 2$).  
 We have
$$ 
0 < \frac{x- y \theta_1}{x- y \theta_2}= \left(\frac{t-\theta_1}{t-\theta_2}\right)^n \left(\frac{\theta_1}{\theta_2} \right)^{-m},
$$
and hence $n$ and $m$ have the same parity. Arguing as for type II solutions, we may write
$$
\log \left( \frac{x- y \theta_1}{x- y \theta_2} \right) = \frac{1}{t^3} + \frac{5}{2t^6} +\frac{25}{3 t^9} + \frac{\kappa_{14}}{t^{12}}, \mbox{ with } \kappa_{14} \in (25.9,26.6),
$$
$$
\; \log \left( \frac{t-\theta_1}{\theta_2-t} \right) = 6 \log t - \frac{\kappa_{15}}{t^3}, \; \mbox{ with } \kappa_{15} \in (2.9,3.1),
$$
$$
\log \left( \frac{\theta_2}{|\theta_1|} \right) = 6 \log t - \frac{2}{t^3} - \frac{\kappa_{16}}{t^6}, \mbox{ for } \kappa_{16} \in (-0.1,0.1),
$$
and consider the identity
$$
\log \left( \frac{x- y \theta_1}{x- y \theta_2} \right) =n  \log \left( \frac{t-\theta_1}{\theta_2-t} \right)
+ m \log \left( \frac{\theta_2}{|\theta_1|} \right).
$$
We thus have
\begin{equation} \label{nm}
n+m = \frac{1}{6 \log t} \left( \frac{2m+\kappa_{15}n+1}{t^3} + \frac{5/2+\kappa_{16} m}{t^6} +\frac{25}{3 t^9} + \frac{\kappa_{14}}{t^{12}} \right).
\end{equation}
It is easy to show from this that $nm \neq 0$ and, indeed, that $n$ and $m$ have opposite signs.
From (\ref{good}), since necessarily $|x- \theta_2 y|<1$, we have $m> 0$ and $n < 0$, say $n=-n_0$.

If $m \geq n_0+2$, then (\ref{nm}) and the inequality $t \geq 10$ imply that
$$
m-n_0 < \frac{1}{6 \log t} \left( \frac{2m- 2.9 n_0 +1 }{t^3}+ \frac{2.6+ 0.1 m}{t^6} \right) < \frac{m-n_0}{2 \,   t^3 \, \log t },
$$
a contradiction. Since $m$ and $n$ have the same parity, we may thus conclude that either $m=n_0$ or $n_0 \geq m+2$. In the first case, we have
$$
\frac{(\kappa_{15}-2) m-1}{t^3} - \frac{5/2+\kappa_{16} m}{t^6} -\frac{25}{3 t^9} - \frac{\kappa_{14}}{t^{12}} =0,
$$
and hence, after a little work, that $m=1, n=-1$, corresponding to the solution $(x,y)=(t^4-2t,1).$
We may therefore assume that  $n_0 \geq  m+2 $ and hence, from (\ref{nm}),
$$
2 \leq n_0 - m < \frac{3.1 n_0- 1.9 m}{6 t^3 \log t} =  \frac{3.1 (n_0-  m)}{6 t^3 \log t} +  \frac{m}{5 t^3 \log t} 
$$
and so
 \begin{equation} \label{rat}
|n| >  m > 4.9 \,  (n_0-m) \, t^3 \log t  \geq 9.8 \, t^3 \log t.
 \end{equation}
 
We may thus conclude, in all cases (i.e. for solutions of type I, II or III), that
\begin{equation} \label{space}
\max \left\{ |m|, |n| \right\} \geq 3.5 \, t^3 \log t.
\end{equation}

\section{ Linear forms in logarithms}

With a lower bound upon $\max \left\{ |m|, |n| \right\}$ in hand, we now turn our attention to
extracting  bounds for expressions of the shape (\ref{rip}).
Our starting point is Siegel's identity :
$$
(\theta_2-\theta_3)(x- y \theta_1)+(\theta_3-\theta_1)(x- y \theta_2)+(\theta_1-\theta_2)(x- y \theta_3)=0.
$$
From this, we have, for example, that 
$ \frac{x-y\theta_3}{x-y\theta_2}$
 not only satisfies equation (\ref{good}), but also
$$
 \frac{x-y\theta_3}{x-y\theta_2} = \frac{\theta_1-\theta_3}{\theta_1-\theta_2} 
 + \frac{\theta_3-\theta_2}{\theta_1-\theta_2} (x - y \theta_1).
 $$
 It follows, if $(x,y)$ is a type I solution, say, that $\frac{x-y\theta_3}{x-y\theta_2}$ and $\frac{\theta_1-\theta_3}{\theta_1-\theta_2}$ are extremely close together, whereby, from (\ref{good}), the linear form
\begin{equation} \label{lam1}
\Lambda_1=\log \left(  \frac{\theta_1-\theta_3}{\theta_1-\theta_2} \right) +n \log \left( \frac{t-\theta_2}{t-\theta_3} \right)-m \log \left( \frac{\theta_2}{\theta_3} \right)
\end{equation}
is necessarily small. Explicitly, we may write $\Lambda_1=\log \left(1+ \tau_1 \right),$ where 
$$
\tau_1=\frac{(\theta_2-\theta_3)(x- y \theta_1)}{(\theta_1-\theta_2)(x- y \theta_3)}.
$$
We have, since $m$ and $n$ are negative in (\ref{unit}) for solutions of type I, with $|m| > 3 |n|$, 
$$
\left \lvert \tau_1\right \rvert=\frac{\theta_3-\theta_2}{\theta_2-\theta_1 } \left(\frac{t-\theta_1}{\theta_3-t}\right)^{n} \left( \frac{\theta_3}{| \theta_1 |} \right)^m<(t^3-3)\left(t^3-3 \right)^{-n} \left( t^9-4 t^6 \right)^{m}  < t^{7.7 m},
$$ 
where, in the last inequality, we are assuming that $(m,n) \neq (-4,-1)$ (a case we treated earlier).
From the fact that $| \log (1+z)| < 2 |z|$, valid for $|z|< 1/2$, we may thus conclude that
 \begin{equation} \label{lowerbound3-1}
  \log \left \lvert \Lambda_1 \right \rvert < \log 2 - 7.7 \, |m| \,  \log t. 
  \end{equation} 

We argue similarly in the case of solutions of types II and III, considering the linear forms
\begin{equation} \label{lam2}
\Lambda_2=\log \left(\frac{\theta_3-\theta_2}{\theta_2-\theta_1} \right)+n \log \left\lvert\frac{t-\theta_1}{t-\theta_3}\right\rvert+ m \log \left\lvert\frac{\theta_3}{\theta_1}\right\rvert
\end{equation}
and
\begin{equation} \label{lam3}
\Lambda_3 = \log \left( \frac{\theta_3-\theta_2}{\theta_3-\theta_1} \right) +n \log \left\lvert\frac{t-\theta_1}{t-\theta_2}\right\rvert+m \log \left\lvert\frac{\theta_2}{\theta_1}\right\rvert.
\end{equation}
Corresponding to (\ref{lowerbound3-1}), we have, after some work, the inequalities
 \begin{equation} \label{flip}
\log  \left| \Lambda_2 \right| < \log 2 -7.9 \, n \log t
 \end{equation}
and
\begin{equation} \label{flip2}
\log  \left| \Lambda_3 \right| < \log 2  -8.9 \, |n| \log t
 \end{equation}

   To find lower bounds for our $|\Lambda_i |$, we require estimates for linear forms in logarithms of algebraic numbers.
 The following is the main result (Theorem 2.1) of Matveev \cite{Mat} :
\begin{proposition} \label{Mat}
Let $\mathbb{K}$ be an algebraic number field of degree $D$ over $\mathbb{Q}$ and set $\chi=1$ if $\mathbb{K}$ is real and $\chi=2$ otherwise. Suppose that $\alpha_1, \alpha_2, \ldots, \alpha_n \in \mathbb{K}^*$ with absolute logarithmic heights $h(\alpha_i)$ for $1 \leq i \leq n$, and that
$$
A_i \geq \max \{ D \, h (\alpha_i), \left| \log \alpha_i \right| \}, \; 1 \leq i \leq n,
$$
for some fixed choice of the logarithm. Define
$$
\Lambda = b_1 \log \alpha_1 + \cdots + b_n \log \alpha_n,
$$
where the $b_i$ are rational integers and set
$$
B = \max \{ 1, \max \{ |b_i| A_i/A_n \; : \; 1 \leq i \leq n \} \}.
$$
Define, further, 
$\Omega =A_1 \cdots A_n$, 
$$
C = \frac{16}{n! \chi} e^n (2n+1+2 \chi) (n+2)(4n+4)^{n+1} \left( e n/2 \right)^{\chi},
$$
$$
C_0 = \log \left( \exp(4.4 n+7) n^{5.5} D^2 \log ( e D) \right)
$$
and
$W_0 = \log \left(
1.5 e B D \log ( e D) \right).$
Then, if $\log \alpha_1, \ldots, \log \alpha_n$ are linearly independent over $\mathbb{Z}$ and $b_n \neq 0$, we have
$$
\log \left| \Lambda \right| > - C \, C_0 \, W_0 \, D^2 \, \Omega.
$$
\end{proposition}

 We can apply this result, with suitable parameter choices, to find  lower bounds upon $\Lambda_i$, for each $i \in \{ 1, 2 , 3 \}$. We will focus our attention on the case of $\Lambda_2$ (where the resulting upper bound upon $t$ is largest). The other cases proceed similarly; details are available from the authors on request.
 
To treat $\Lambda_2$, we choose
  $$
  \alpha_3= \frac{\theta_3-\theta_2}{\theta_2-\theta_1}, \; \alpha_2 =\left\lvert\frac{t-\theta_1}{t-\theta_3}\right\rvert , \; \alpha_1 =\left\lvert\frac{\theta_3}{\theta_1}\right\rvert, \;b_3=1, \; b_2=n, \; b_1=m \mbox{ and } D=6.
  $$
From (\ref{estimation}), we have
$$
h \left(\frac{\theta_3-\theta_2}{\theta_3-\theta_1}\right) \leq 2h \left(\theta_3-\theta_2 \right) =\frac{2}{3} \log \left((\theta_3-\theta_2)(\theta_3-\theta_1)(\theta_2-\theta_1)\right) < 6 \log t ,
$$
$$
h \left(\frac{\theta_3}{\theta_1}\right) =\frac{1}{6} \log \left(\frac{\theta_3}{\theta_1}\right)^2 <3 \log t
$$
and
$$
h \left(\frac{t-\theta_1}{t-\theta_3}\right) =\frac{1}{6} \log \left(\frac{t-\theta_3}{t-\theta_2}\right) ^2<3\log t.
$$
Therefore we can take 
$$
A_3=36 \log t, \;  A_1=A_2= 18 \log t \; \mbox{ and }\; B=\frac{n}{2}
$$
in Proposition \ref{Mat} to conclude that
$$
 \log \, \lvert \Lambda_2 \rvert >-8.4 \cdot 10^{15} \log^3 t  \; \log \,(35 n).
 $$
 Combining this with (\ref{flip}), we thus have that
 \begin{equation} \label{lug}
 \frac{n}{\log\, ( 35 n)} < 1.07 \cdot 10^{15} \log^2 t.
 \end{equation}
Appealing to (\ref{lower}),  we may therefore conclude that $t \leq 576241$ (so that $n < 8.9 \cdot 10^{18}$). 


\subsubsection{Small values of t}

At this point, there are a number of ways to proceed to handle the remaining values of $t$.
One of which, which would be particularly valuable if our bound upon $t$ was less good, would be to observe that we can rewrite our linear form $\Lambda_2$ as
$$
\Lambda_2 = m \log \alpha_2 - \log \alpha_1,
$$
where 
$$
\alpha_1=\left (\frac {\theta_2-\theta_1} {\theta_3-\theta_2} \right ) \left\lvert\frac {t-\theta_3}{t-\theta_1}\right\rvert^{k+1}
\; \mbox{ and } \; \alpha_2=\left\lvert\frac{t-\theta_1}{t-\theta_3}\right\rvert^3 \left\lvert\frac{\theta_3}{\theta_1}\right\rvert.
$$
Since (\ref{help}) ensures that $k$ is small, the height of $\alpha_1$ is not too large and so we can profitably apply lower bounds for linear forms in two logarithms, rather than three; typically, this leads to much improved numerical results.
    
In our case, since the values of $t$ under consideration are not especially large, we will instead appeal to a result from Diophantine approximation. Specifically, we will use a lemma of Mignotte \cite{Mig}, a variant of a classical result by Baker and Davenport \cite{BaDa}.

\begin{lemma} \label{baker-Davenport}
Let $\Lambda= \mu \alpha + \nu \beta +\delta $ where $\alpha, \beta$ and $\delta$ are nonzero real numbers and where $\mu, \nu $ are rational integers with $\left\lvert \mu \right\rvert < A.$
Let $Q>0$ be a real number and suppose that $\gamma_1$ and $\gamma_2$ satisfy 
$$
\left\lvert \gamma_1-\frac{\alpha}{\beta}\right\rvert <\frac{1}{100 \, Q^2} \quad \text{and} \quad  \left\lvert \gamma_2-\frac{\delta}{\beta}\right\rvert <\frac{1}{Q^2}.
$$ 
 Further, let $p/q$ be a rational number with $1 \leq q \leq Q $ and $\left\lvert \gamma_1-p/q \right\rvert <q^{-2}$, and suppose $q \left\lvert\left\lvert q \gamma_2 \right\rvert\right\rvert \geq 1.01A+2$, where $\left\lvert\left\lvert \cdot \right\rvert\right\rvert $ denotes the distance to nearest integer. Then
  \begin{equation} \label{mignotteinq}
   \left\lvert\Lambda \right\rvert > \frac{\left\lvert \beta \right\rvert}{Q^2}.
   \end{equation}
  \end{lemma}
   To apply this lemma in our situation, we choose
  $$
  \alpha= \log \left\lvert\frac{\theta_3}{\theta_1}\right\rvert, \;  \beta=\log \left\lvert\frac{t-\theta_1}{t-\theta_3}\right\rvert, \; \delta=\log \left(  \frac{\theta_3-\theta_2}{\theta_2-\theta_1} \right),  \; \mu=m, \; 
  \nu=n
  $$  
  and $A = 3 \cdot 10^{18}$.  For each $t$ with $10 \leq t \leq 576241$, we compute $\alpha, \beta$ and $\delta$ with suitable precision, choosing $Q=10^{60}$. In each case, searching the continued fraction expansion of $\gamma_1$, we are able to find a convergent  $p/q$ with the desired properties. We may thus conclude that $|\Lambda_2| > 10^{-120}$. On the other hand, combining (\ref{lower}) and (\ref{flip}), we have
  $$
  \log | \Lambda_2 | < -27.65 \, t^3 \log^2 t,
  $$
  an immediate contradiction.  
Full details of this computation are available from the authors upon request.
      
\section{Closing remarks} 

It is a relatively routine exercise nowadays to solve a given cubic Thue equation (or even a fairly large number of them).  Algorithmic routines, based on a paper of Tzanakis and de Weger \cite{TdW}, building on the foundational work of Baker \cite{Bak}, exist in a number of computer algebra packages, including PARI and Magma. We used code based on the latter to solve the equations  corresponding to (\ref{main-eq}) with $|t| < 10$ (where our previous arguments may fail), as well as all cubic Thue equations of the shape $F(x,y)=1$ for forms of positive discriminant $D_F \leq 10^7$. For the latter, we argued as in \cite{Ben2}, appealing to work of Belabas \cite{Be} and Belabas and Cohen \cite{BC}, \cite{BC2}.

For such  forms, there are precisely $9$ equivalence classes with $N_F \geq 5$ which are inequivalent to $F_{i,t} (x,y)$ for any $i \in \{ 1, 2, 3 \}$ and $t \in \mathbb{Z}$. These are

$$
\begin{array}{|ccl|} \hline
F(x,y) & D_F & N_F  \\ \hline
x^3-3xy^2+y^3 & 81 & 6 \\
x^3 + x^2 - 3 x - 1 & 148 & 5 \\
x^3+2x^2y-5xy^2+y^3 & 361 & 6  \\
x^3 - 5 x - 1 & 473 & 5 \\
x^3 - 7 x - 1 & 1345 & 5 \\
x^3 + 9 x^2 - 12 x - 21 & 108729 & 5 \\
x^3 + 21 x^2 - 2 x - 21 & 783689 & 5 \\
x^3 + 21 x^2 - x - 22 & 810661 & 5 \\
x^3 + 18 x^2 - 21 x - 37 & 1257849 & 5 \\
 \hline
\end{array}
$$
\vskip1ex

It is also worth noting that the families of forms $F_{1,t_1}, F_{2,t_2}$ and $F_{3,t_3}$ are essentially disjoint. It is a routine exercise (since the discriminants in each family are essentially squares, whereby we may apply Runge's method) to show that the only cases where we have 
$$
F_{i,t_i} (x,y) \sim F_{j,t_j} (x,y) \; \mbox{ with } \; i, j \in \{ 1, 2, 3 \}, \; i \neq j
$$
correspond to
$$
F_{1,0}(x,y) \sim F_{1,2}(x,y) \sim F_{2,1} (x,y) \sim F_{3,1}(x,y), \; \; F_{1,4}(x,y) \sim F_{3,-1}(x,y)
$$
and 
$$
F_{2,0}(x,y) \sim F_{3,0}(x,y).
$$
Here,  since $F_{1,-t}(x,y)=F_{1,t-1}(y,x)$ and $F_{2,t}(x,y)=F_{2,-t}(x,y)$ we suppose that $t_1, t_2 \geq 0$.

The solution of the remaining cases of the equation $F_{2,t}(x,y)=1$, i.e those with $0 \leq t < 1.35 \cdot 10^{14}$, is within computational range nowadays, requiring slight refinements of the arguments of \cite{Wa5}, in conjunction with appeal to state-of-the-art lower bounds for linear forms in three complex logarithms, due to Mignotte \cite{Migster}. The reason why this family is somewhat harder to solve is that the lower bounds upon the growth of the exponents of the corresponding fundamental units is less strong in this case. Indeed one has only that the larger exponent is of size $\gg t^{3/2}$. To obtain this, Wakabayashi \cite{Wa5} employs local arguments reminiscent of Skolem's $p$-adic method. 


\end{document}